\input amstex
\documentstyle{amsppt}

\magnification=\magstep1

\define\e{\varepsilon}

\define\th{\theta}

\define\conv{\text{conv}}
\define\etc{, \dots ,}
\define \sumi{\sum_{i=1}^M}

\define\rn{$\Bbb R^n \,$}

\define\Rn{\Bbb R^n}
\define\bn{$B^n_2 \,$}
\define\Bn{B^n_2}
\define\nor #1{\left \| #1 \right \|}
\define\enor #1{\Bbb E \, \nor{#1}}
\define\tens #1{#1 \otimes #1}
\define\pr#1#2{\langle {#1} , {#2} \rangle}
\define\ba#1#2{\Cal B_{#2} ( {#1} )}

\topmatter

\title
Random vectors in the isotropic position
\endtitle

\author
M. Rudelson
\endauthor

\thanks
Research at MSRI is supported in part  by NSF grant DMS-9022140.
\endthanks

\affil
MSRI \ Texas A \& M University
\endaffil

\address
Mathematical Sciences Research Institute, \hfill \break
1000 Centennial Drive, Berkeley, CA 94720, USA
\endaddress

\email
mark\@msri.org
\endemail

\vskip 3cm
\abstract
Let $y$ be a random vector in \rn, satisfying 
$$
\Bbb E \, \tens{y} = id.
$$
Let $M$ be a natural number and let $y_1 \etc y_M$ be independent copies
 of $y$. 
We prove that for some absolute constant $C$
$$
\enor{\frac{1}{M} \sumi \tens{y_i} - id} \le 
C \cdot \frac{\sqrt{\log M}}{\sqrt{M}} \cdot 
\left ( \enor{y}^{\log M} \right )^{1/ \log M},
$$
provided that the last expression is smaller than 1.

We apply this estimate to obtain a new proof of a result of Bourgain 
concerning the number of random points needed to bring a convex body
into a nearly isotropic position.
\endabstract

\endtopmatter

\document

\head 1. Introduction \endhead
The problem we consider has arisen from a question studied by 
R.~Kannan, L.~Lov\'asz and M.~Simonovits \cite{K-L-S}.
To construct a fast algorithm for calculating the volume of a convex body,
 they needed to bring it into some symmetric position.
More precisely, let $K$ be a convex body in \rn. 
We shall say that it is in the isotropic position if for any $x \in \Rn$
$$
\frac{1}{\text{vol }(K)} \int_K \pr{x}{y}^2 \, dy = \nor{x}^2
$$
By $\nor{\cdot}$ we denote the standard Euclidean norm.

The notion of isotropic position was extensively studied by V.~Milman and
A.~Pajor \cite{M-P}.
Note that our definition is consistent with \cite{K-L-S}.
The normalization in \cite{M-P} is slightly different.

If the information about the body $K$ is uncomplete it is impossible to
bring it exactly to the isotropic position.
So, the definition of the isotropic position has to be modified to allow 
a small error.
We shall say that the body $K$ is in $\e$-isotropic position if
for any $x \in \Rn$
$$
(1-\e) \cdot \nor{x}^2 \le 
\frac{1}{\text{vol }(K)} \int_K \pr{x}{y}^2 \,   dy \le (1+\e) \cdot \nor{x}^2.
$$
Let $\e>0$ be given.
Consider $M$ random points $y_1 \etc y_M$ independently uniformly distributed 
in $K$ and put
$$
T=\frac{1}{M}  \sumi \tens{y_i}.
$$
If $M$ is sufficiently large, than with high probability 
$$
\nor{T- \frac{1}{\text{vol} (K)} \int_K \tens{y}}
$$
will be small, so the body $T^{-1/2} K$ will be in $\e$-isotropic position.
R.~Kannan, L.~Lov\'asz and M.~Simonovits (\cite{K-L-S}) proved that it is
enough to take 
$$
M=c \frac{n^2}{\e}
$$
for some absolute constant $c$.
This estimate was significantly improved by J.~Bourgain \cite{B}. 
Using rather delicate geometric considerations he has shown that one can take 
$$
M =C(\e) n \, \log^3 n.
$$
Since the situation is invariant under a linear transformation, we may assume
that the body $K$ is in the isotropic position.
Then the result of Bourgain may be reformulated as follows:
\proclaim{Theorem 0}\cite{B}
Let $K$ be a convex body in \rn in the isotropic position.
Fix $\e>0$ and choose independently $M$ random points $x_1 \etc x_M \in K$,
$$
M \ge  C(\e) n \, \log^3 n.
$$
Then with probability at least $1-\e$ for any $x \in \Rn$ one has
$$
(1-\e)\nor{x}^2 \le \frac{1}{M} \sumi \pr{x}{y}^2 \le (1+\e)\nor{x}^2.
$$
\endproclaim
We shall show that this theorem follows from a general result about random 
vectors in \rn.
Let $y$ be a random vector.
Denote by $\Bbb E \, X$ the expectation of a random variable $X$. 
We say that $y$ is in the isotropic position if
$$
\Bbb E \, \tens{y} = id. \tag 1.1
$$
If $y$ is uniformly distributed in a convex body $K$, then this is equivalent
to the fact that $K$ is in the isotropic position.

We prove the following
\proclaim{Theorem 1}
Let $y \in \Rn$ be a random vector in the isotropic position. 
Let $M$ be a natural number and let $y_1 \etc y_M$ be independent 
copies of $y$.
Then
$$
\enor{\frac{1}{M} \sumi \tens{y_i} - id} \le 
C \cdot \frac{\sqrt{\log M}}{\sqrt{M}} \cdot 
\left ( \enor{y}^{\log M} \right )^{1/ \log M}, \tag 1.2
$$
provided that the last expression is smaller than 1.
\endproclaim
Here and later $C,c$, etc. denote absolute constants whose values may vary
from line to line.

\demo{Remark}
Taking the trace of (1.1) we obtain that $\enor{y}^2=n$, so  to make the right
hand side of (1.2) smaller than 1, we have to assume that $M \ge cn \log n$.
\enddemo
Using Theorem 1 we prove a better estimate of the length of approximate 
John's decompositions \cite{R1} and thus improve the results about
approximating a convex body by another one having a small number of contact 
points, obtained in \cite{R2}.
Estimating the moment of the norm of random vector in a convex body,
we obtain a different proof of Theorem 0 which gives also a better estimate.

\head 2. Main results. \endhead
The proof of Theorem 1 consists of two steps.
First we introduce a Bernoulli random process  and estimate the expectation 
of the norm in (1.2) by the expectation of its supremum.
Then we construct a majorizing measure to obtain a bound for the latest.

The first step is relatively standard.
Let $\e_1 \etc \e_M$ be independent Bernoulli variables  taking values 
$1,-1$ with probability $1/2$ and let 
$y_1 \etc y_M, \quad \bar y_1 \etc \bar y_M$ be independent copies of $y$.
Denote $\Bbb E_y, \, \Bbb E_{\e}$ the expectation according to $y$ and $\e$
respectively.
Since $\tens{y_i}-\tens{\bar y_i}$ is a symmetric random variable, we have
$$
\align
&\Bbb E_y \, \nor{\frac{1}{M} \sumi \tens{y_i} - id} \le
\Bbb E_y \, \Bbb E_{\bar y} \, \nor{\frac{1}{M} \sumi \tens{y_i} - 
\frac{1}{M} \sumi \tens{\bar y_i}} = \\
&\Bbb E_{\e} \,\Bbb E_y \,\Bbb E_{\bar y} \, 
\nor{\frac{1}{M} \sumi \e_i ( \tens{y_i}- \tens{\bar y_i})} \le
2 \, \Bbb E_y \, \Bbb E_{\e} \, \nor{\frac{1}{M} \sumi \e_i \tens{y_i}}.
\endalign
$$
To estimate the last expectation, we need the following Lemma, which 
generalizes Lemma 1 \cite{R3}.
\proclaim{Lemma}
Let $y_1 \etc y_M$ be vectors in \rn and let $\e_1 \etc \e_M$ be 
independent Berno\-ulli variables taking values $1,-1$ with probability $1/2$.
Then 
$$
\enor{\sumi \e_i \tens{y_i}} \le C \sqrt{\log M} \cdot 
\max_{i=1 \etc M} \nor{y_i} \cdot \nor{\sumi \tens{y_i}}^{1/2}.
$$
\endproclaim
We postpone the proof of the Lemma to the next section.

Applying the Lemma, we get
$$
\aligned
&\enor{\frac{1}{M} \sumi \tens{y_i} - id} \le \\
C \cdot \frac{\sqrt{\log M}}{M} \cdot 
&\left ( \Bbb E \, \max_{i=1 \etc M}  \nor{y_i}^2 \right )^{1/2} \cdot 
\left ( \enor{\sumi \tens{y_i}}\right )^{1/2}. 
\endaligned  \tag 2.1
$$
We have 
$$
\align
\left ( \Bbb E \, \max_{i=1 \etc M}  \nor{y_i}^2 \right )^{1/2} \le
&\left ( \Bbb E \, 
\left ( \sumi \nor{y_i}^{\log M} \right )^{2/ \log M} \right )^{1/2} 
\le \\
&M^{1/ \log M} \cdot  \left (\enor{y}^{\log M} \right )^{1/ \log M} .
\endalign
$$
Thus, denoting
$$
D=\enor{\frac{1}{M} \sumi \tens{y_i} - id},
$$
we obtain by (2.1)
$$
D \le C \cdot\frac{\sqrt{\log M}}{\sqrt{M}} \cdot 
\left (\enor{y}^{\log M} \right )^{1/ \log M} \cdot (D+1)^{1/2}.
$$
If 
$$
C \cdot\frac{\sqrt{\log M}}{\sqrt{M}} \cdot 
\left (\enor{y}^{\log M} \right )^{1/ \log M} \le 1,
$$
we get
$$
D \le 2C \cdot\frac{\sqrt{\log M}}{\sqrt{M}} \cdot 
\left (\enor{y}^{\log M} \right )^{1/ \log M},
$$
which completes the proof of Theorem 1.  
\medpagebreak

We turn now to the applications of Theorem 1.
Applying Theorem 1 to the question of Kannan, Lov\'asz and Simonovits, 
we obtain the following 
\proclaim{Corollary 2.1}
Let $\e >0$ and let $K$ be an $n$-dimensional convex body in the isotropic 
position.
Let 
$$
M \ge C \cdot \frac{n}{\e^2} \cdot \log^2 \frac{n}{\e^2}
$$
and let $y_1 \etc y_M$ be independent random vectors uniformly distributed 
in $K$.
Then
$$
\enor{\frac{1}{M} \sumi \tens{y_i} - id} \le \e.
$$
\endproclaim
\demo{Proof}
It follows from  a result of S.~Alesker \cite{A}, that 
$$
\Bbb E \, \exp \left ( \frac{\nor{y}^2}{c \cdot n} \right ) \le 2
$$
for some absolute constant $c$.
Then
$$
\align
\enor{y}^{\log M} \le
&\left (\Bbb E \,  e^{ \frac{\nor{y}^2}{c \cdot n} } \right )^{1/2} \cdot
\left ( \Bbb E \, \left ( \nor{y}^{2 \log M} \cdot
e^{- \frac{\nor{y}^2}{c \cdot n} } \right ) \right )^{1/2}\le \\
&\sqrt{2} \cdot 
\left ( \max_{t \ge 0} t^{\log M} \cdot 
e^{-\frac{t}{c \cdot n}} \right )^{1/2} \le
(C \cdot n \cdot \log M )^{\frac{\log M}{2}}. 
\endalign
$$
Corollary 2.1 follows from this estimate and Theorem 1. \qquad \qed
\enddemo

By a Lemma of Borell \cite{M-S, Appendix III}, most of the volume of a
convex body in the isotropic position is concentrated within the Euclidean 
ball of radius $c \sqrt{n}$.
So, it might be of interest to consider a random vector uniformly 
distributed in the intersection of a convex body and such a ball.
In this case the previous estimate may be improved as follows.
\proclaim{Corollary 2.2}
Let $\e,R >0$ and let $K$ be an $n$-dimensional convex body in the isotropic 
position.
Suppose that $R \ge c \sqrt{\log 1/\e}$ and let 
$$
M \ge C_0 \cdot \frac{R^2 \cdot n}{\e^2} \cdot \log \frac{R^2 \cdot n}{\e^2}
\tag 2.2
$$
and let $y_1 \etc y_M$ be independent random vectors uniformly distributed 
in  \break
${K \cap R \sqrt{n}) \cdot \Bn}$.
Then
$$
\enor{\frac{1}{M} \sumi \tens{y_i} - id} \le \e.
$$
\endproclaim

\demo{Proof}
Denote $a= R \cdot \sqrt{n}$ and
let $z$ be a random vector uniformly distributed in $K \cap a \Bn$.
Then for $x\in \Bn$
$$
\align
\Bbb E \, \pr{z}{x}^2 = \frac{\text{vol }(K)}{\text{vol }(K \cap a \Bn)} \cdot
\Big ( \frac{1}{\text{vol }(K)} &\int_K \pr{y}{x}^2 \, dy \ -  \\
 \frac{1}{\text{vol }(K)} &\int_K \pr{y}{x}^2 \cdot 
\bold 1_{\{u \bigm | \nor{u} \ge a  \}}(y) \, dy  \Big ).
\endalign
$$
By a result of S.~Alesker \cite{A} and Khinchine type inequality \cite{M-P},
we have
$$
\frac{\text{vol }(K)}{\text{vol }(K \cap a \Bn)} \le 1+e^{-ca^2/n} \le
1+ \frac{\e}{4}
$$
and
$$
\align
&\frac{1}{\text{vol }(K)} \int_K \pr{y}{x}^2 \cdot 
\bold 1_{\{u \bigm | \nor{u} \ge a  \}}(y) \, dy \le \\
&\left ( \frac{1}{\text{vol }(K)} \int_K \pr{y}{x}^4 \, dy \right)^{1/2} \cdot
\left ( \frac{1}{\text{vol }(K)} 
\int_K\bold 1_{\{u \bigm | \nor{u} 
\ge a \}}(y) \, dy \right )^{1/2} \le \\
&C e^{-ca^2/2n} \le \frac{\e}{4}.
\endalign
$$
Thus for any $x \in \Bn$
$$
| \Bbb E \,\pr{z}{x}^2 -1 | \le \frac{\e}{2}.
$$
Define a random vector 
$$
y= ( \Bbb E \, \tens{z} )^{-1/2} z.
$$
Then $y$ is in the isotropic position and
$$
\left ( \enor{y}^{\log M} \right )^{1/ \log M} \le
\nor{( \Bbb E \, \tens{z} )^{-1/2}} \cdot 
\left ( \enor{z}^{\log M} \right )^{1/ \log M} \le 2a,
$$
so
$$
\enor{\frac{1}{M} \sumi \tens{y_i} - id} \le 
C \cdot \frac{\sqrt{\log M}}{\sqrt{M}} \cdot 2a \le \frac{\e}{2}
$$
provided the constant $C_0$ in (2.2) is large enough.
Thus,
$$
\enor{\frac{1}{M} \sumi \tens{z_i} - id} \le 
\enor{\frac{1}{M} \sumi \tens{y_i} - id} \cdot \nor{\Bbb E \, \tens{z}} +
\nor{\Bbb E \, \tens{z} -id} \le \e. 
$$
\qed
\enddemo

The next application is connected to the approximation of a convex body by 
another one having a small number of contact points \cite{R2}.
Let $K$ be a convex body in \rn such that the ellipsoid of minimal 
volume containing it is the standard Euclidean ball \bn.
Then by the theorem of John, there exist 
$N \le (n+3)n/2$   points $z_1 ,\dots, z_N \in K, \quad \nor{x_i} =1$ and 
$N$ positive numbers $c_1, \dots, c_N$ satisfying the following 
system of equations
$$
\align
id &= \sum_{i=1}^N c_i \, z_i \otimes z_i \tag 2.3 \\
0 &= \sum_{i=1}^N c_i \, z_i.              \tag 2.4
\endalign
$$
It was shown in \cite{R1} for convex symmetric bodies and in \cite{R2}
in the general case, that  the identity operator can be approximated by a sum
of a smaller number of terms $\tens{x_i}$.
We derive from Theorem 1 the following corollary, which improves Lemma 3.1 
\cite{R2}.
\proclaim{Corollary 2.3}
Let $\e>0$ and let $K$ be a convex body in \rn, so that the ellipsoid of 
minimal volume containing it is \bn.
Then there exist
$$
M \le \frac{C}{\e^2} \cdot n \cdot \log \frac{n}{\e} \tag 2.5
$$
contact points $x_1, \dots, x_M$ and a vector 
$u, \ \nor{u} \le \frac{C}{\sqrt{M}}$, 
so that the identity operator in \rn has the following representation
$$
\align
id=\frac{n}{M} \, & \sumi  \tens{(x_i+u)}+S,  \\
\intertext{where} 
&\sumi (x_i+u) = 0 \tag 2.6\\
\intertext{and} 
&\nor{S : \ell^n_2 \to \ell^n_2} < \e. 
\endalign
$$
\endproclaim
\demo{Proof}
Let (2.3) be a decomposition of the identity operator.
Let $y$ be a random vector in \rn, taking values $\sqrt{n} z_i$ with 
probability $c_i/ \sqrt{n}$.
Then, by (2.3), $y$ is in the isotropic position.
Obviously, for all $1 \le p < \infty$ 
$$
\left ( \enor{y}^p \right )^{1/ p} = \sqrt{n}.
$$
So, taking $M$ as in (2.5), we obtain that for sufficiently large $C$
$$
\nor{\frac{1}{M} \sumi \tens{y_i} - id} \le \frac{\e}{2} \tag 2.7
$$
with probability greater than $3/4$.
Since by (2.4), $\Bbb E \, y=0$ and $\nor{y}=\sqrt{n}$, we have
$$
\nor{\sumi y_i} \le 2 \sqrt{M} \tag 2.8
$$
with probability greater than $3/4$.
Take $y_1 \etc y_M$ for which (2.7) and (2.8) hold and put 
$$
x_i= \frac{1}{\sqrt{n}} \cdot y_i, \qquad \qquad u = -\frac{1}{M} \sumi x_i.
$$
Then (2.6) is satisfied and 
$$
\align
&\nor{\frac{n}{M} \,  \sumi  \tens{(x_i+u)} -id} \le \\
&\nor{\frac{n}{M} \, \sumi \tens{x_i} - id} + n \cdot \nor{\tens{u}} \le 
\frac{\e}{2} + \frac{4\, n}{M} \le \e. \qed
\endalign
$$  
\enddemo

Substituting Lemma 3.1 \cite{R2} by Corollary 2.3 in the proof of Theorem 1.1 
\cite{R2} we obtain the following
\proclaim{Corollary 2.3}
 Let B be a convex body in \rn and let $\e>0$.
There exists a convex body $K \subset \Rn$, so that ${\text d}(K,B) \le 1+\e$
and the number of contact points of $K$ with the ellipsoid of minimal volume 
containing it  is less than
$$ 
M(n, \e)= \frac{C}{\e^2} \cdot n \cdot \log \frac{n}{\e}.
$$
\endproclaim

\head 3. Proof of the Lemma \endhead
The proof of the Lemma is similar to that of Lemma 1 \cite{R3}.
For the reader's convenience we present here a complete proof.

Without loss of generality, we may assume that 
$$
\nor{\sumi \tens{y_i}} =1.
$$
Define a random process
$$
V_x= \sumi \e_i \pr{x}{y}^2
$$
for $x \in \Bn$.
We have to estimate
$$
\Bbb E \, \sup_{x \in \Bn} V_x. \tag 3.1
$$
Note that the process $V_x$ has a subgaussian tail estimate
$$
\Cal P \{ V_x - V_{\bar x} > a \} \le 
\exp \left ( -C \frac{a^2}{\tilde d^2 (x, \bar x)} \right ),
$$
where $C$ is an absolute constant and 
$$
\tilde d (x, \bar x) =
\left ( \sumi \Big ( \pr{x}{y_i}^2 -\pr{\bar x}{y_i}^2 \Big )^2 \right )^{1/2}.
$$
The function $\tilde d$ is not a metric on \bn, since 
$\tilde d (x, \bar x) =0$ does not imply $x=\bar x$. 
To avoid this obstacle, we shall estimate $\tilde d$ by a quasimetric $d$
defined by
$$
d (x, \bar x) = \left ( \sumi \pr{x-\bar x}{y_i}^2 
\Big ( \pr{x}{y_i}^2 + \pr{\bar x}{y_i}^2 \Big ) \right )^{1/2}.
$$
Then for all $, \bar x \in \Bn$
$$
\tilde d (x, \bar x) \le \sqrt{2} \cdot d (x, \bar x),
$$
so we may treat $V_x$ as a subgaussian process with the quasimetric $d$.
It can be easily shown that $d$ satisfies a generalized triangle inequality
$$
d (x, \bar x) \le 4 \cdot ( d (x, z) +  d (z, \bar x) ) \tag 3.2
$$
for all $x,\bar x, z \in \Bn$.

Denote by $\ba{x}{\rho}$ a ball in the quasimetric $d$ with center $x$
and radius $\rho$.
Then for any $x \in \Bn$ and $\rho>0$ we have 
$$
\text{conv } \ba{x}{\rho} \subset \ba{x}{4 \rho}.
$$
The proof of this fact is the same as that of Lemma 3 \cite{R3}, so we shall
omit it.
To estimate the expression (3.1), we apply the following version of the 
Majorizing measure theorem.
\proclaim{Theorem}
Let $(T,d)$ be a quasimetric space.
Let $(X_t )_{t \in T}$ be a collection of mean 0 random variables with the
subgaussian tail estimate
$$
 \Cal P \, \{| X_t - X_{\bar t}| > a \} \le 
\exp \left ( -c \frac{a^2}{d^2 (t, \bar t)} \right ),
$$
for all $a>0$.
Let $r>1$ and let $k_0$ be a natural number so that the diameter of $T$ is 
less than $r^{-k_0}$.
Let $\{\varphi_k \}_{k=k_0}^{\infty}$ be a sequence of functions from 
$T$ to $\Bbb R^+$,  uniformly bounded by a constant depending only on $r$.
Assume that there exists $\sigma >0$ so that for any $k$ the functions 
$\varphi_k$ satisfy the following condition: \hfil \break
\par\flushpar
 for any $s \in T$ and for any points 
$t_1 \etc t_N \in B_{r^{-k}}(s)$ with mutual distances at least $r^{-k-1}$ 
one has
$$
\max_{j=1 \etc N} \varphi_{k+2}(t_j) \ge \varphi_k(s) + \sigma \cdot r^{-k}
\cdot \sqrt{\log N}.  \tag 3.3
$$
Then 
$$
\Bbb E \sup_{t \in T} X_t \le C(r) \cdot \sigma^{-1}.
$$
\endproclaim
This Theorem is a combination of the majorizing
measure theorem of Fernique \cite{L-T} and the general majorizing
measure construction of Talagrand (Theorems 2.1 and 2.2 \cite{T1} or
Theorems 4.2, 4.3 and Proposition 4.4 \cite{T2}).

Let 
$$
Q= \max_{i= 1 \etc M} \nor{y_i}.
$$
Let $r$ be a natural number and let $k_0$ and $k_1$ be the largest numbers 
so that
$$
\align
&r^{-k_0} \ge Q \\
&r^{-k_1} \ge \frac{Q}{\sqrt{n}}.
\endalign
$$

Then $k_1-k_0 \le (2 \log r )^{-1} \log n$.
Define now the functions $\varphi_k: \Bn \to \Bbb R$ by
$$
\alignat 2
\varphi_k(x)=& \min \{ \nor{u}^2 \ \Big | 
\ u \in \conv \ba{x}{8 r^{-k}} \} +
\frac{k-k_0}{\log M}, 
&& \qquad \text{if } k=k_0 \etc k_1, \\
\varphi_{k}(w)=& 1+ \frac{1}{2 \log r} + \sum_{l=k_1}^k r^{-l} \cdot 
\frac{\sqrt{n \cdot \log ( 1+ 4 Q r^l  )}}
{Q \cdot \sqrt{\log M}},
&& \qquad \text{if } k>k_1.
\endalignat
$$
The functions $\varphi_k$ form a nonnegative nondecreasing sequence bounded
by a constant depending only on $r$.
Indeed, for $k \le k_1$, 
$$
\varphi_k(x) \le 1+ \frac{1}{2 \log r}\cdot \frac{\log n}{\log M}.
$$
For $k>k_1$ we have 
$$
\align
\varphi_k(w) \le & 1+ \frac{1}{2 \log r} +\sum_{l=k_1}^{\infty} r^{-l} \cdot 
\frac{\sqrt{n \cdot \log ( 1+4 Q r^l )}}{Q \sqrt{\log M}}  \le \\
& 1+ \frac{1}{2 \log r} + c(r) \cdot r^{-k_1} \cdot \sqrt{n} \cdot
\frac{\sqrt{\log ( 1+4 Q r^{k_1} )}} {Q \sqrt{\log M}} \le C(r).
\endalign
$$
Now let $x_1 \etc x_N \in \ba{x}{r^-k}$ and suppose that
$$
d(x_i, x_j) > r^{-k-1}
$$
for any $i \neq j$.
We have to prove that
$$
\max_{j=1 \etc N} \varphi_{k+2} (x_i) \ge \varphi_k (x) + 
\frac{C}{Q \sqrt{\log M}} \cdot \sqrt{\log N}. \tag 3.4
$$
Note that for $x, \bar x \in \Bn$
$$
\align
d(x, \bar x) \le &\sqrt{2} \cdot \max_{i=1 \etc M} | \pr{x- \bar x}{y_i} |
\cdot \left ( \sumi \pr{x}{y_i}^2 + \pr{\bar x}{y_i}^2 \right )^{1/2} \le \\
&\sqrt{2} \cdot \max_{i=1 \etc M} | \pr{x- \bar x}{y_i} | \cdot
\nor{\sumi \tens{y_i}} \cdot \left ( \nor{x}^2 + \nor{\bar x}^2 \right )^{1/2} 
\le \\
&2 \cdot \max_{i=1 \etc M} | \pr{x- \bar x}{y_i} |.
\endalign
$$
Define a norm in \rn by 
$$
\nor{x}_Y = \max_{i=1 \etc M} | \pr{x}{y_i} |.
$$
If $k \ge k_1-2$ then (3.4) follows from a simple entropy estimate.
Indeed, we have 
$$
\align
N \le &N(\Bn,d,r^{-k-1}) \le N(\Bn, \nor{\cdot}_Y, \frac{1}{2}r^{-k-1}) \le \\
&N(\Bn, \nor{\cdot}, \frac{1}{2 Q}r^{-k-1}) \le (1+ 4 Q \cdot r^{k+1})^n,
\endalign
$$
since $\nor{\cdot}_Y \le Q \nor{\cdot}$.
Suppose now that $k<k_1 -2$.
For $j=1 \etc N$ denote $z_j$ the point of $\conv \ba{x_j}{8r^{-k-2}}$
for which the minimum of $\nor{z}^2$ is attained and denote
$u$ the similar point of $\ba{x}{8r^{-k}}$.
Put
$$
\th =\max_{j=1\etc N} \nor{z_j}^2 - \nor{u}^2.
$$
We have to show that 
$$
r^{-k} \cdot \left ( c \cdot Q \cdot \sqrt{\log M} \right )^{-1} \cdot
\sqrt{\log N}
\le \max_{j=1 \etc N} \varphi_{k+2}(x_j)-\varphi_k(x)= \th+ \frac{2}{\log M}.
\tag 3.5
$$

Since $d(x_i, x_j) \ge r^{-k-1}$,  it follows from (3.2) that
$$
d(z_i, z_j) \ge \frac{1}{2} r^{-k-1},
$$
provided $r$ is sufficiently large.
From the other side,
$$
d(x, z_j) \le 4 \Big (  d(x,x_j) + d(x_j,z_j) \Big ) \le 8 r^{-k}.
$$
Since $\dfrac{z_j + u}{2} \in \conv \ba{x}{8r^{-k}}$, 
and $\nor{u} \le \nor{z_j}$, we have
$$
\nor{\frac{z_j-u}{2}}^2 = \frac12 \nor{z_j}^2 + \frac12 \nor{u}^2 -
\nor{\frac{z_j+u}{2}}^2 \le \nor{z_j}^2 - \nor{\frac{z_j+u}{2}}^2 \le
\nor{z_j}^2 - \nor{u}^2,
$$
so,
$$
\nor{z_j-u} \le 2 \sqrt{\th}. \tag 3.6
$$
Thus, $N$ is bounded by the $\frac12 r^{-k-1}$-entropy of the set 
$K=u+2\sqrt{\th} \Bn$ in the quasimetric $d$.
To estimate this entropy we partition the set $K$ into $S$ disjoint subsets
having diameter less than $\delta =\frac{1}{16} r^{-k-1} \th^{-1/2}$ in the 
$\nor{\cdot}_Y$ metric.

Let $g$ be a Gaussian vector in \rn, normalized by $\enor{g}^2=n$.
Denote by $N(B,\Delta, \e)$ $\e$-entropy of the set $B$ in the metric $\Delta$.
By dual Sudakov minoration \cite{L-T} we have 
$$
\align
\sqrt{\log S} \le &\sqrt{\log N(2 \sqrt{\th} \Bn, \nor{\cdot}_Y,\delta)} \le \\
&\frac{c}{\delta} \cdot 2 \sqrt{\th} \cdot \enor{g}_Y \le
C \cdot r^k \th \cdot \Bbb E \max_{i=1 \etc M} | \pr{g}{y_i} | \le
C \cdot r^k \th \cdot Q \cdot \sqrt{\log M}.  \tag 3.7
\endalign 
$$

If $S \ge \sqrt{N}$, we are done, because in this case (3.7) implies (3.5).
Suppose that  $S \le \sqrt{N}$.
Then there exists an element of the partition containing at least $\sqrt{N}$
points $z_j$.
Let $J \subset \{ 1 \etc N \}$ be the set of the indices of these points.
We have 
$$
\nor{z_j-z_l}_Y \le \frac{1}{16} r^{-k-1} \cdot \th^{-1/2}  \tag 3.8
$$
for all $j,l \in J$.

For $j=1 \etc M$ denote 
$$
I_j= \{ i \in \{ 1 \etc M \} \bigm | |\pr{z_j}{y_i} \ge 2 \, | \pr{u}{y_i}| \}.
$$
Then (3.6) imlies that
$$
\sum_{i \in I_j} \pr{z_j}{y_i}^2 \le 
2 \sum_{i \in I_j} \pr{z_j-u}{y_i}^2 + 2 \sum_{i \in I_j} \pr{u}{y_i}^2 \le
8 \th + \frac{1}{2} \sum_{i \in I_j} \pr{z_j}{y_i}^2,
$$
so,
$$
\sum_{i \in I_j} \pr{z_j}{y_i}^2 \le 16 \th. \tag 3.9
$$
Since $d(z_j,z_l) \ge \frac12 r^{-k-1}$, we have
$$
\align
 \left ( \frac12 r^{-k-1} \right )^2 \le &\sumi \pr{z_j-z_l}{y_i}^2
\cdot \Big ( \pr{z_j}{y_i}^2+\pr{z_l}{y_i}^2 \Big ) \le \\
& \sumi  \pr{z_j-z_l}{y_i}^2 \cdot 4 \pr{u}{y_i}^2 \ + \\
&\max_{i \in I_j} \pr{z_j-z_l}{y_i}^2 \cdot
\sum_{i \in I_j} \pr{z_j}{y_i}^2 + 
\max_{i \in I_l} \pr{z_j-z_l}{y_i}^2 \cdot
\sum_{i \in I_l} \pr{z_l}{y_i}^2.
\endalign
$$

Combining (3.8) and (3.9) we get that the last expression is bounded by
$$
2 \cdot 16 \th \cdot \left ( \frac{\th^{-1/2}}{8} r^{-k-1} \right )^2 +
4 \sumi \pr{z_j-z_l}{y_i}^2 \cdot \pr{u}{y_i}^2.
$$
Define a norm $\nor{\cdot}_{\Cal E}$ by
$$
\nor{x}_{\Cal E}= 
\left ( \sumi \pr{x}{y_i}^2 \cdot \pr{u}{y_i}^2 \right )^{1/2}.
$$
Then, for all $j,l  \in J, \ j \ne l$ we have 
$$
\nor{z_j-z_l}_{\Cal E} \ge \frac18 r^{-k-1}.
$$
Applying again dual Sudakov minoration, we obtain
$$
\align
\sqrt{\log |J|} \le 
&\sqrt{\log N(2 \sqrt{\th} \Bn, \nor{\cdot}_{\Cal E}, \frac18 r^{-k-1})} \le
c r^k \cdot 2\sqrt{\th} \cdot \enor{g}_{\Cal E} \le \\
&c r^k \cdot 2\sqrt{\th} \cdot 
\left ( \Bbb E \, \sumi \pr{g}{y_i}^2 \cdot \pr{u}{y_i}^2 \right )^{1/2} \le \\
&C r^k \cdot 2\sqrt{\th} \cdot \max_{i=1 \etc M} \nor{y_i} \cdot 
\nor{\sumi \tens{y_i}} \cdot \nor{u} \le C r^k \cdot 2\sqrt{\th} \cdot Q.
\endalign
$$
Since for all $\th>0$
$$
2 \sqrt{\th} \le \sqrt{\log M} \cdot \th +\frac{1}{\sqrt{\log M}},
$$
we get
$$
\sqrt{\log N} \le 2 \sqrt{\log |J|} \le  C  \cdot Q \cdot r^k
\cdot \sqrt{\log M} \cdot \left ( \th +\frac{1}{\log M} \right ), 
$$
so (3.5) is satisfied.

\enddocument